\documentclass[12pt,reqno]{amsart}

\def\divides{{\mathchoice{\mathrel{\bigm|}}{\mathrel{\bigm|}}{\mathrel{|}}{\mathrel{|}}}}
\def\notdivides{\mathrel{\kern-3pt\not\!\kern3.5pt\bigm|}}
\def\smallnotdivides{\mathrel{\kern-2pt\not\!\kern3.5pt\vert}}
\def\Divides{\divides\!\divides}

\renewcommand{\mid}{\colon}
\newtheorem{theorem}{Theorem}
\newtheorem{proposition}[theorem]{Proposition}
\newtheorem{lemma}[theorem]{Lemma}

\theoremstyle{definition}

\newtheorem{example}[theorem]{Example}

\theoremstyle{remark}

\newcommand{\ord}{\operatorname{ord}}
\newcommand{\fix}{\operatorname{\mathcal F}}
\newcommand{\least}{\operatorname{\mathcal L}}
\newcommand{\orbit}{\operatorname{\mathcal O}}

\begin{document}

\bibliographystyle{amsplain}

\title[Orbit counting with
an isometric direction]{Orbit counting with\\an isometric direction}

\author{G. Everest}
\author{V. Stangoe}
\author{T. Ward}

\address{School of Mathematics, University of East
  Anglia, Norwich NR4 7TJ, United Kingdom}
\email{t.ward@uea.ac.uk}

\thanks{The second author thanks the EPSRC for
their support via a Doctoral Training Award, and
the third author
thanks Richard Sharp for telling
him about~\cite{MR98h:53055}, and the
Max Planck Institute in Bonn where this
paper was completed.}

\renewcommand{\subjclassname}{\textup{2000} Mathematics Subject Classification}
\begin{abstract}
Analogues of the prime number theorem and Merten's theorem are
well-known for dynamical systems with hyperbolic behaviour. In
this paper a $3$-adic extension of the circle doubling map is
studied. The map has a $3$-adic eigendirection in which it behaves
like an isometry, and the loss of hyperbolicity leads to weaker
asymptotic results on orbit counting than those obtained for
hyperbolic maps.
\end{abstract}

\maketitle

\section{Introduction}

As part of the analogy between the behaviour of closed orbits in
hyperbolic dynamical systems and prime numbers
in~\cite{MR92f:58141} both the prime number
theorem~\cite{MR85c:58089}, \cite{MR85i:58105} and Merten's
theorem~\cite{MR93a:58142} have corresponding results. Roughly
speaking, these results mean the following for a hyperbolic
transformation $f:X\to X$ with topological entropy $h=h(f)$. A
{\sl closed orbit} $\tau$ of {\sl length} $\vert\tau\vert=n$ is a
set of the form
$$
\{x,f(x),f^2(x),\dots,f^n(x)=x\}
$$
with cardinality $n$. The analogue of the prime number theorem
states that
\begin{equation}\label{hyppnt}
\pi_f(X)=\#\{\tau\mid\vert\tau\vert\le X\}\sim
\frac{e^{h(X+1)}}{X(e^h-1)},
\end{equation}
and the analogue of Merten's theorem states that
\begin{equation}\label{hypmerten}
\sum_{\vert\tau\vert\le X}\frac{1}{e^{h\vert\tau\vert}}\sim
\log X+C.
\end{equation}
The results in the papers cited are more precise, and apply to
flows as well as hyperbolic maps. We will use the word
`hyperbolicity' rather loosely, and in particular will apply it to
group endomorphisms whose natural invertible extension is an
expansive automorphism of a finite-dimensional group. If one
allows $p$-adic as well as complex eigenvalues, this corresponds
to having no eigenvalues of unit size, and taking the invertible
extension does not affect the number of periodic orbits.

Without hyperbolicity, less is known. For ergodic toral
automorphisms that are not hyperbolic (the {\sl quasihyperbolic}
automorphisms of~\cite{MR84g:28017}),
Waddington~\cite{MR92k:58219} has shown an analogue
of~\eqref{hyppnt} and Noorani~\cite{MR2001k:37036} an analogue
of~\eqref{hypmerten}. In both cases additional terms appear, but
the basic shape remains the same, and the technique of exploiting
a meromorphic extension of the dynamical zeta function beyond its
radius of convergence still works. In a different direction,
Knieper~\cite{MR98h:53055} finds asymptotic upper and lower bounds
for the function counting closed geodesics on rank-$1$ manifolds
of non-positive curvature of the form
\begin{equation*}
A\frac{e^{hX}}{X}
\leq\#\{\mbox{closed geodesics of length} \le X\}\leq
B e^{hX}
\end{equation*}
for constants $A,B > 0$.

Our purpose here is to describe for a specific map
how~\eqref{hyppnt} and~\eqref{hypmerten} change when hyperbolicity
is lost via the introduction of an isometric extension that gives
the local structure of the dynamical system a (non-Archimedean)
isometric direction. This work is a special case of more general
results in~\cite{stangoethesis}, where many isometric directions,
a positive-characteristic analogue, and more refined estimates are
considered. The essential issues that arise are well illustrated
by the map studied here, which avoids many technicalities. In
particular, the bad behaviour of the dynamical zeta function
appears for this simple example already.

\section{Orbits and periodic points}

We first recall some basic properties of maps and their periodic
points (see~\cite{MR2002i:11026} for example). Let $T:X\to X$ be a
map. The number of points of period $n$ under $T$ is
\begin{equation*}
\fix_n(T)=\#\{x\in X\mid T^nx=x\},
\end{equation*}
the number of periodic points with least period $n$ under $T$ is
\begin{equation*}
\least_n(T)=\#\{x\in X\mid T^nx=x\mbox{ and
}\#\{T^kx\}_{k\in\mathbb N}=n\},
\end{equation*}
and the number of closed orbits of length $n$ is
\begin{equation}\label{orbitleast}
\orbit_n(T)=\least_n(T)/n.
\end{equation}
Since
\begin{equation}\label{sumdiv}
\fix_n(T)=\sum_{d\vert n}\least_d(T),
\end{equation}
the sequence of values taken by any one of
these quantities determines the others.
A consequence of this is a relationship
between the number of orbits of a map and the
number of orbits of its iterates which will be
needed in a special case later.

\begin{lemma}\label{You were only waiting for this moment to arise}
For any $k\ge1$,
$
\orbit_n(T^k)=\frac{1}{n}\sum_{d\vert n}\mu({n}/{d})
\sum_{d^{\prime}\vert dk}d^{\prime}\orbit_{d^{\prime}}(T).
$
In particular,
\begin{equation*}
\orbit_n(T^2)=\begin{cases}
2\orbit_{2n}(T)+\orbit_{n}(T)&\mbox{if }n\mbox{ is odd,}\\
2\orbit_{2n}(T)&\mbox{if }n\mbox{ is even}.
\end{cases}
\end{equation*}
\end{lemma}

\begin{proof}
Clearly $\fix_d(T^k)=\fix_{dk}(T)$ for all $d\ge1$
so, by M{\"o}bius inversion of~\eqref{sumdiv},
\begin{eqnarray*}
\least_n(T^k)&=&\sum_{d\vert n}\mu(n/d)\fix_d(T^k)
=\sum_{d\vert n}\mu(n/d)\fix_{kd}(T)\\
&=&\sum_{d\vert n}\mu(n/d)\sum_{d^{\prime}\vert dk}
\least_{d^{\prime}}(T),
\end{eqnarray*}
which gives the result by~\eqref{orbitleast}.

The expression for $\orbit_n(T^2)$ may be checked
as follows. First notice that by~\eqref{sumdiv}
it is sufficient to check that the given expression
gives the right value to $\fix_n(T^2)$.
If $n$ is odd, then all factors of $n$ are odd so
\begin{eqnarray*}
\fix_n(T^2)&=&\sum_{d\divides n}d\orbit_d(T^2)\\
&=&\sum_{d\divides n}2d\orbit_{2d}(T)+\sum_{d\divides n}d\orbit_d(T)\\
&=&\sum_{d\divides 2n}d\orbit_d(T)=\fix_{2n}(T)
\end{eqnarray*}
as required.
If $n$ is even, then
\begin{eqnarray*}
\fix_n(T^2)&=&\sum_{d\divides n}d\orbit_d(T^2)\\
&=&\sum_{2\divides d\divides n}2d\orbit_{2d}(T)
+\sum_{2\smallnotdivides d\divides n}2d\orbit_{2d}(T)
+\sum_{2\smallnotdivides d\divides n}d\orbit_d(T)\\
&=&\sum_{d\vert n}2d\orbit_{2d}(T)+
\sum_{2\smallnotdivides d\divides n}d\orbit_d(T)\\
&=&\sum_{2\divides d\divides 2n}d\orbit_d(T)+
\sum_{2\smallnotdivides d\divides 2n}d\orbit_d(T)\\
&=&\sum_{d\divides 2n}d\orbit_d(T)=\fix_{2n}(T).
\end{eqnarray*}
\end{proof}

\section{A $\mathbb Z_3$-extension of the circle-doubling map}

Consider the map $\phi:x\mapsto 2x$ on the ring $\mathbb
Z[\frac{1}{3}]$. Write $X=\widehat{\mathbb Z[\frac{1}{3}]}$ for
the dual (character) group, and $f=\widehat{\phi}$ for the dual
map. The pair $(X,f)$ is the map we will study here, using the
following properties from~\cite{MR99b:11089}
and~\cite{MR99k:58152}.
\begin{itemize}
\item $X$ is a compact metrizable group and $f$ is a continuous
endomorphism of $X$. \item Locally, the action of $f$ is isometric
to the map $(s,t)\mapsto(2s,2t)$ on an open set in $\mathbb
R\times\mathbb Q_3$. \item In this local picture, the real
coordinate is stretched by $2$ while the action on the $3$-adic
coordinate is an isometry. \item The number of points of period
$n$ under $f$ is $\left\vert2^n-1\right\vert\cdot
\left\vert2^n-1\right\vert_3\!.$ \item The topological entropy of
$f$ is $\log 2$.
\end{itemize}
The map $f$ may also be thought of as an extension of the familiar
circle-doubling map by a cocycle taking values in $\mathbb Z_3$,
the $3$-adic integers. The short exact
sequence\begin{equation}\label{shortexact}\textstyle0\longrightarrow\mathbb
Z\overset{\imath}\longrightarrow\mathbb
Z[\frac{1}{3}]\longrightarrow\mathbb Z[\frac{1}{3}]/\mathbb
Z\longrightarrow0\end{equation}(where $\imath$ is the inclusion
map) commutes with $x\mapsto 2x$. The dual of~\eqref{shortexact}
is the short exact sequence$$0\longrightarrow\mathbb
Z_3\longrightarrow X\overset{\hat{\imath}}\longrightarrow\mathbb
T\longrightarrow0$$which commutes with the dual of $x\mapsto 2x$.
This sequence expresses $f:X\to X$ as an extension of the
circle-doubling map by a cocycle taking values in $\mathbb Z_3$.
This extension kills certain periodic orbits.

\section{Dynamical zeta function}

From now on $f:X\to X$ is the map from the previous section.
The arguments below are elementary: the most sophisticated
number-theoretic fact needed is that the divisors of an
odd number are all odd.

First notice that
\begin{eqnarray*}
\vert2^n-1\vert_3&=&\vert
(3-1)^n-1\vert_3\\
&=&\vert3^n-n3^{n-1}+\dots+(-1)^{n-1}3n+(-1)^n-1\vert_3\\
&=&\begin{cases}
\frac{1}{3}\vert n\vert_3&\mbox{if }n\mbox{ is even,}\\
1&\mbox{if }n\mbox{ is odd.}
\end{cases}
\end{eqnarray*}
In
particular,
\begin{equation}\label{blackbird}
\vert4^n-1\vert_3=\vert2^{2n}-1\vert_3=\textstyle\frac{1}{3}\vert2n\vert_3=\frac{1}{3}\vert
n\vert_3.
\end{equation}
Since $\vert n\vert_3=
3^{-\ord_3(n)}\ge3^{-\log_3(n)}\ge1/n$, it follows that
\begin{equation}\label{polybound}
\textstyle\frac{1}{3n}\le\vert2^n-1\vert_3\le 1\mbox{ for all
}n\ge1,
\end{equation}
so just as for the circle doubling map the logarithmic growth rate
of periodic points gives the topological entropy,
$$
\frac{1}{n}\log\fix_n(f)\longrightarrow\log2=h(f).
$$
Thus the radius of convergence of the
dynamical zeta function
$$
\zeta(z)=\exp\sum_{n=1}^{\infty}\frac{z^n}{n}\fix_n(f)
$$
is $\exp\left(-h(f)\right)=\frac{1}{2}$.

The bound~\eqref{polybound} means that the number of periodic
points for the circle doubling map is only polynomially larger
than the number of periodic points for $f$. This is enough room to
make a real difference: for example, in contrast to the hyperbolic
case, $\frac{\fix_{n+1}(f)}{\fix_n(f)}$ does not converge as
$n\to\infty$ (indeed, for $S$-integer systems in zero
characteristic, this convergence essentially characterizes
hyperbolicity, cf.~\cite[Th.~6.3]{MR99b:11089}).

One of the key tools in the hyperbolic case is a meromorphic
extension of the zeta function to a larger disc; in our setting
this is impossible.

\begin{proposition} The dynamical zeta function
of $f$ has natural boundary $\vert z\vert=\frac{1}{2}$.
\end{proposition}

\begin{proof}
Let $\xi(z)=\sum_{n=1}^{\infty}\frac{z^n}{n}
\vert2^n-1\vert\cdot
\vert2^n-1\vert_3$ so
$\zeta(z)=\exp(\xi(z))$.
Now
\begin{eqnarray*}
\xi(z)&=&\sum_{n=0}^{\infty}
\frac{z^{2n+1}}{2n+1}
(2^{2n+1}-1)+\sum_{n=1}^{\infty}
\frac{z^{2n}}{2n}(2^{2n}-1)\vert2^{2n}-1\vert_3\\
&=&\log\left(\frac{1-z}{1-2z}\right)
-{\frac{1}{2}}\log\left(
\frac{1-z^2}{1-4z^2}\right)+
\sum_{n=1}^{\infty}
\frac{z^{2n}}{2n}(2^{2n}-1)\vert2^{2n}-1\vert_3.
\end{eqnarray*}
Write $\frac{1}{6}\xi_1(z)$ for the last term in this
expression, so
\begin{eqnarray*}
\xi_1(z)&=&3
\sum_{n=1}^{\infty}
\frac{z^{2n}}{n}(4^{n}-1)\vert4^{n}-1\vert_3\\
&=&
\sum_{n=1}^{\infty}
\frac{z^{2n}}{n}(4^n-1)\vert n\vert_3
\end{eqnarray*}
by~\eqref{blackbird}. We shall show that $\xi_1(z)$ has infinitely many
logarithmic singularities on the circle $|z|=\frac{1}{2}$,
each of which corresponds to a zero of $\zeta(z)$.

Write $3^a\Divides n$ to mean that $3^a\divides n$ but
$3^{a+1}\notdivides n$. Notice that $3^a\Divides n$ if and only if
$\vert n\vert_3=3^{-a}$. Then $\xi_1$ may be split up according to
the size of $\vert n\vert_3$ as
\begin{eqnarray*}
\xi_1(z)&=&\sum_{j=0}^{\infty}\frac{1}{3^j}\sum_{3^j\Vert n}
\frac{z^{2n}}{n}(4^n-1)\\
&=&\sum_{j=0}^{\infty}\frac{1}{3^j}\eta_j^{(4)}(z),
\end{eqnarray*}
where
$$
\eta_j^{(a)}(z)=\sum_{3^j\Divides n}
\frac{z^{2n}}{n}(a^n-1).
$$
Then
\begin{eqnarray*}
\eta_0^{(a)}(z)&=&\sum_{3^0\Divides n}\frac{z^{2n}}{n}(a^n-1)\\
&=&\sum_{n=1}^{\infty}\frac{z^{2n}}{n}(a^n-1)
-\sum_{n=1}^{\infty}\frac{z^{6n}}{3n}(a^{3n}-1)\\
&=&\log\left(\frac{1-z^2}{1-az^2}\right)-
\frac{1}{3}\log\left(\frac{1-z^6}{1-a^3z^6}\right){\!\!},
\end{eqnarray*}
\begin{equation*}
\eta_1^{(4)}(z)=\sum_{3^1\Divides n}\frac{z^{2n}}{n}(4^n-1)
=\sum_{3^0\Divides n}\frac{z^{6n}}{3n}(4^{3n}-1)
=\textstyle\frac{1}{3}\eta_0^{(4^3)}(z^3),
\end{equation*}
\begin{equation*}
\eta_2^{(4)}(z)=\textstyle\frac{1}{9}\eta_0^{(4^9)}(z^9),
\end{equation*}
and so on.
Thus
\begin{equation*}
\xi_1(z)=\log\left(\frac{1-z^2}{1-(2z)^2}\right)+
2\sum_{j=1}^{\infty}\frac{1}{9^j}\log\left(\frac{1-(2z)^{2\times3^j}}{1-z^{2\times3^j}}
\right){\!\!},
\end{equation*}
so
\begin{equation*}
\left\vert\zeta(z)\right\vert=
\left\vert\frac{1-z}{1-2z}\right\vert\cdot
\left\vert\frac{1-(2z)^{2\vphantom{3^j}}}{1-z^2}\right\vert^{1/2}\!\!\!\!\cdot
\left\vert\frac{1-z^{2\vphantom{3^j}}}{1-(2z)^{2}}\right\vert^{1/6}\!\!\!\!\cdot
\prod_{j=1}^{\infty}\left\vert
\frac{1-(2z)^{2\times3^j}}{1-z^{2\times3^j}}\right\vert^{1/3\times9^j}
\end{equation*}

It follows that the series defining $\zeta(z)$ has a zero at all
points of the form $\frac{1}{2}e^{2\pi ij/3^r}$, $r\ge1$ so $\vert
z\vert=\frac{1}{2}$ is a natural boundary for $\zeta(z)$.
\end{proof}

A natural boundary appears for most of the
$S$-integer dynamical systems studied in~\cite{MR99b:11089},
but only {\it ad hoc} proofs of this exist.
Natural boundaries also arise for certain
`random' zeta functions~\cite{MR2003g:37034} and
in zeta functions for higher-rank actions~\cite{MR97e:58185}.
The appearance of singularities at all $g^r$th unit roots
creating a natural boundary has been exploited by
Hecke and Mahler in other contexts~\cite{hecke}, \cite{MR89f:30006}.

\section{Prime orbit theorem}

The first result shows the asymptotics
for $\pi_f(X)$, and is arrived at by comparison
with the well-understood asymptotics for
the circle-doubling map
$g:\mathbb T\to\mathbb T$,
$x\mapsto2x$ mod $1$. Here
$\fix_n(g)=2^n-1$,
$\orbit_n(g)=\frac{1}{n}\sum_{d\vert n}
\mu(n/d)(2^d-1)$, and~\eqref{hyppnt} applies
to show that
\begin{equation}\label{pntstated}
\pi_g(X)=\sum_{n\le X}\frac{1}{n}\sum_{d\vert n}
\mu(n/d)(2^d-1)\sim\frac{2^{X+1}}{X}.
\end{equation}
Pursuing the analogy between results like~\eqref{pntstated}
and the prime number theorem, the next result is
analogous to Tchebycheff's theorem.

\begin{theorem}\label{All your life}
\begin{equation*}
\pi_f(X)\le\pi_g(X)
\mbox{ for all }X\ge 1,
\end{equation*}
and
\begin{equation*}
\limsup_{X\to\infty}\frac{X\pi_f(X)}{2^{X+1}}\le1,
\quad
\liminf_{X\to\infty}\frac{X\pi_f(X)}{2^{X+1}}\ge\frac{1}{3}.
\end{equation*}
\end{theorem}

\begin{proof} Let $b_n=2^n-1$ and
$a_n=b_n\vert 2^n-1\vert_3$, so
$$
\pi_f(X)=\sum_{n\le X}\orbit_n(f)=
\sum_{n\le X}\frac{1}{n}\sum_{d\vert n}\mu(n/d)a_d.
$$
We first claim that
\begin{equation}\label{Take these broken wings and learn to fly}
\orbit_n(f)\le\orbit_n(g)\mbox{ for all }n\ge1.
\end{equation}
Notice that this does not follow {\it a priori}
from the fact that
$$
\fix_n(f)\le\fix_n(g)\mbox{ for all }n\ge1;
$$
it is easy to construct pairs of maps with one of these
inequalities but not the other (see Example~\ref{now I'm
freefalling}). However,~\eqref{Take these broken wings and learn
to fly} does hold here because of the exponential growth in $a_n$.

For odd $n$, $a_n=b_n$ so $\orbit_n(f)=\orbit_n(g)$ (since
all factors of $n$ are also odd).

Now assume that $n$ is even and notice that
\begin{equation}\label{Blackbird singing in the dead of night}
\sum_{d\vert n,d<n}(2^d-1)\le\textstyle\frac{2}{3}(2^n-1)\mbox{ for all }
n\ge1.
\end{equation}
It follows that
\begin{eqnarray*}
\least_n(f)&=&\sum_{d\vert n}\mu(n/d)a_d\\
&\le&a_n+\sum_{d\vert n,d<n}a_d\\
&\le&\textstyle\frac{1}{3}b_n+\sum_{d\vert n,d<n}b_d\\
&\le&\sum_{d\vert n}\mu(n/d)b_d=\least_n(g)
\end{eqnarray*}
by~\eqref{Blackbird singing in the dead of night}.
Thus $\orbit_n(f)\le\orbit_n(g)$ for all $n\ge1$,
so~\eqref{Take these broken wings and learn to fly}
-- and hence
the upper bound in Theorem~\ref{All your life} -- is proved.

Turning to the lower bound, write
$$
\delta(X)=\pi_g(X)-\pi_f(X)\ge0,
$$
and notice that
\begin{equation*}
\delta(X)=\sum_{n\le X}\left(\orbit_n(g)-\orbit_n(f)\right)
=\sum_{2\vert n\le X}\left(\orbit_n(g)-\orbit_n(f)\right)
\le
\sum_{2\vert n\le X}\orbit_n(g).
\end{equation*}
So we need to estimate the size of
$\sum_{2\vert n\le X}\orbit_n(g)$.
Notice that
$$\fix_n(g^2)=4^n-1$$
and $g^2$ (the map $x\mapsto4x$ mod 1 on the
circle) is hyperbolic, so
$$
\sum_{n\le X}\orbit_n(g^2)\sim\frac{4^{X+1}}{3X}.
$$
By Lemma~\ref{You were only waiting for this moment to arise},
\begin{equation*}
\orbit_n(g^2)=\begin{cases}
2\orbit_{2n}(g)+\orbit_{n}(g)&\mbox{if }n\mbox{ is odd,}\\
2\orbit_{2n}(g)&\mbox{if }n\mbox{ is even}.
\end{cases}
\end{equation*}
So
\begin{eqnarray*}
2\sum_{n\le X}\orbit_{2n}(g)&=&
\sum_{2\smallnotdivides n\le X}\left(
\orbit_n(g^2)-\orbit_n(g)\right)+
\sum_{2\divides n\le X}\orbit_n(g^2)\\
&=&\sum_{n\le X}\orbit_n(g^2)-
\sum_{2\smallnotdivides n\le X}\orbit_n(g).
\end{eqnarray*}
Now $\sum_{n\le X}\orbit_n(g^2)\sim\frac{4^{X+1}}{3X}$,
so
$$
\sum_{n\le X}\orbit_{2n}(g)\sim\frac{2}{3}\cdot
\frac{4^X}{X}-\frac{1}{2}\sum_{2\smallnotdivides n\le X}
\orbit_n(g).
$$
On the other hand, $\sum_{n\le X}\orbit_n(g)\sim\frac{2^{X+1}}{X}$,
so the last term is of lower order. It follows that
$$
\sum_{n\le X}\orbit_{2n}(g)\sim\frac{2}{3}\cdot\frac{4^X}{X},
$$
so
$$
\pi_g(X)-\pi_f(X)\le
\sum_{2\divides n\le X}\orbit_{n}(g)\sim
\frac{2}{3}\cdot\frac{4^{X/2}}{X/2}=\frac{2}{3}\cdot
\frac{2^{X+1}}{X}.
$$
Thus
\begin{equation*}
\limsup_{X\to\infty}\frac{X\pi_f(X)}{2^{X+1}}\le1,
\quad
\liminf_{X\to\infty}\frac{X\pi_f(X)}{2^{X+1}}\ge\frac{1}{3}.
\end{equation*}
\end{proof}

Finding the exact values of the upper and lower limits requires
more careful estimates; numerical evidence suggests the limit does
not exist, and moreover that the sequence
$\left(\frac{X\pi_f(X)}{\vphantom{A^{A^x}}2^{X+1}}\right)_{X\ge1}$
has more than two limit points.

In the proof of Theorem~\ref{All your life} we mentioned that for
a pair of maps $f,g$, $\fix_n(f)\le\fix_n(g)$ for all $n\ge1$ does
not imply that $\orbit_n(f)\le\orbit_n(g)$. Of course
$$\orbit_n(f)\le\orbit_n(g)\mbox{ for all
}n\ge1\implies\fix_n(f)\le\fix_n(g)\mbox{ for all }n\ge1$$
by~\eqref{orbitleast} and~\eqref{sumdiv}.

\begin{example}\label{now I'm freefalling}
Given any sequence $(o_n)_{n\ge1}$ of non-negative integers, there
is a map with $o_n$ closed orbits of length $n$ (indeed, there is
a $C^{\infty}$ map of the 2-torus with this property
by~\cite{alistair}). Let $f$ be chosen with $\orbit_1(f)=1$,
$\orbit_2(f)=3$ and $\orbit_n(f)=0$ for $n\ge3$. Similarly, let
$g$ be chosen with $\orbit_1(g)=6$, $\orbit_2(g)=1$ and
$\orbit_n(g)=0$ for $n\ge3$. Then $\orbit_n(f)$ is not bounded
above by $\orbit_n(g)$, but
$$\fix_n(f)=4+3\cdot(-1)^n<\fix_n(g)=7+(-1)^n\mbox{ for
all }n\ge1.$$
\end{example}

\section{Merten's theorem}

For the circle doubling map $g$, the dynamical analogue of
Merten's theorem says (in a weak formulation)
\begin{equation}\label{You were only waiting for this moment to be free}
\sum_{n\le X}\frac{\orbit_n(g)}{2^n}
=\log X+O(1).
\end{equation}

In our setting, this asymptotic is lost
in much the same way as Theorem~\ref{All your life}
loses the single asymptotic~\eqref{hyppnt}

\begin{theorem}\label{Take these sunken eyes and learn to see}
$$
\textstyle\frac{1}{2}\log{X}+O(1)\le \sum_{n\le
X}\frac{\orbit_n(f)}{2^n} \le\log X+O(1).
$$
\end{theorem}

\begin{proof}
From the proof of Theorem~\ref{All your life},
$\orbit_n(f)\le\orbit_n(g)$ for all
$n\ge1$, so the upper asymptotic
is an immediate consequence
of~\eqref{You were only waiting for this moment to be free}.

To get a lower bound, notice that
$$
\sum_{n\le X}\frac{\orbit_n(f)}{2^n}=\sum_{2\smallnotdivides n\le X}
\frac{\orbit_n(f)}{2^n}+
\sum_{2\divides n\le X}
\frac{\orbit_n(f)}{2^n}\ge
\sum_{2\smallnotdivides n\le X}
\frac{\orbit_n(g)}{2^n}=A(X).
$$
This may be estimated using
Lemma~\ref{You were only waiting for this moment to arise}
again.
First,
\begin{equation*}
A(X)+\sum_{2\divides n\le X}\frac{\orbit_n(g)}{2^n}=
\sum_{n\le X}\frac{\orbit_n(g)}{2^n}=\log X+O(1).
\end{equation*}
By Lemma~\ref{You were only waiting for this moment to arise},
\begin{eqnarray*}
2\sum_{2\divides n\le X}\frac{\orbit_n(g)}{2^n}&=&
2\sum_{m\le X/2}\frac{\orbit_{2m}(g)}{2^{2m}}\\
&=&\sum_{2\smallnotdivides m\le
 X/2}\frac{\orbit_{m}(g^2)-\orbit_{m}(g)}{2^{2m}}
+\sum_{2\divides m\le X/2}
\frac{\orbit_{m}(g^2)}{2^{2m}}\\
&=&\sum_{m\le X/2}\frac{\orbit_m(g^2)}{2^{2m}}-
\sum_{2\smallnotdivides m\le X/2}\frac{\orbit_m(g)}{2^{2m}}\\
&\ge&\log\frac{X}{2}+O(1)=\log X +O(1)
\end{eqnarray*}
since $\orbit_n(g)\le 2^n$ implies that
$$
\sum_{2\smallnotdivides m\le X/2}\frac{\orbit_m(g)}{2^{2m}}\le
\sum_{m\le X/2}\frac{\orbit_m(g)}{2^{2m}}<\infty.
$$
\end{proof}

Just as in Theorem~\ref{All your life}, numerical
evidence suggests that the lower asymptotic is
genuinely lower than $\log X$.


\providecommand{\bysame}{\leavevmode\hbox to3em{\hrulefill}\thinspace}
\providecommand{\MR}{\relax\ifhmode\unskip\space\fi MR }
\providecommand{\MRhref}[2]{%
  \href{http://www.ams.org/mathscinet-getitem?mr=#1}{#2}
}
\providecommand{\href}[2]{#2}

\end{document}